\documentclass[11pt,numbers]{elsarticle}

\usepackage{amsmath,amsfonts,amssymb,amsthm,hyperref,mathtools,graphicx}

\usepackage{enumitem} 
\usepackage{appendix} 
\usepackage{dsfont} 
\usepackage{xcolor} 
\usepackage{multirow} 
\usepackage{multicol} 
\usepackage{booktabs} 
\usepackage{geometry} 
\usepackage{eso-pic}

\newtheorem{theorem}{Theorem}
\newtheorem{proposition}[theorem]{Proposition}
\newtheorem{lemma}[theorem]{Lemma}
\newtheorem{corollary}[theorem]{Corollary}

\theoremstyle{definition}
\newtheorem{remark}{Remark}

\numberwithin{equation}{section}
\numberwithin{theorem}{section}
\numberwithin{table}{section}
\numberwithin{figure}{section}

\newcommand{\N}{\mathbb{N}}
\newcommand{\R}{\mathbb{R}}
\newcommand{\EE}{\mathsf{E}} 
\newcommand{\Bias}{\mathsf{Bias}} 
\newcommand{\Var}{\mathsf{Var}} 
\newcommand{\bb}[1]{\boldsymbol{#1}}
\newcommand{\OO}{\mathcal{O}}
\newcommand{\oo}{\mathrm{o}}
\newcommand{\rd}{\mathrm{d}}
\newcommand{\ind}{\mathds{1}}
\newcommand{\e}{\varepsilon}

\allowdisplaybreaks

\begin{document}

\AddToShipoutPictureBG*{
  \AtPageUpperLeft{
    \put(\LenToUnit{0.5\paperwidth},\LenToUnit{-1.5cm}){ 
      \makebox[0pt][c]{ 
        \fbox{
          \begin{minipage}{50em}
              \small \centering
              This manuscript was accepted for publication in Statistics (Taylor \& Francis).\\
              This version {\it may differ} from the published version (\href{https://doi.org/10.1080/02331888.2026.2698050}{doi:10.1080/02331888.2026.2698050}) in typographic details.
          \end{minipage}
        }
      }
    }
  }
}

\begin{frontmatter}

\title{On the Dirichlet-kernel Gasser--M\"uller estimator and its competitors \\for fixed design regression on the simplex}

\author[a1]{Hanen Daayeb}
\author[a2]{Christian Genest}
\author[a1]{Salah Khardani}
\author[a3]{Nicolas Klutchnikoff}
\author[a4]{Fr\'ed\'eric Ouimet}

\address[a1]{Laboratoire des équations aux dérivées partielles, Universit\'e de Tunis El-Manar, Tunisia}
\address[a2]{Department of Mathematics and Statistics, McGill University, Montr\'eal (QC) Canada}
\address[a3]{Univ Rennes, CNRS, IRMAR - UMR 6625, F-35000 Rennes, France}
\address[a4]{D\'epartement de math\'ematiques et d'informatique, Universit\'e du Qu\'ebec \`a Trois-Rivi\`eres, \newline Trois-Rivi\`eres (QC) Canada}

\cortext[mycorrespondingauthor]{Corresponding author. Email address: Frederic.Ouimet2@uqtr.ca}

\begin{abstract}
A Dirichlet-kernel Gasser--M\"uller (D-GM) estimator is introduced for fixed design regression on the simplex, extending the univariate analog due to Chen [\textit{Statist. Sinica}, vol.~10(1) (2000), pp. 73–91]. Its pointwise bias and variance, asymptotic normality, and mean integrated squared error are investigated. Some simulation experiments are conducted to compare its small-sample performance with that of two recently proposed alternatives: the Dirichlet-kernel Nadaraya--Watson (D-NW) and local linear (D-LL) estimators. The simulation results reveal that the D-LL estimator is best among the D-LL, D-NW, and D-GM estimators and that the proposed D-GM estimator is worst. A real data analysis is also reported for the GEMAS dataset to analyze the relationship between soil composition and pH levels across various agricultural and grazing lands in Europe.
\end{abstract}


\end{frontmatter}

\thispagestyle{empty}

\vspace{-2mm}
\section{Introduction}

Regression analysis is a fundamental tool for investigating how a response variable $Y$ depends on a $d$-dimensional vector of explanatory variables $\bb{X} = (X_1, \dots, X_d)$. In the nonparametric setting, kernel-based methods are frequently used \citep{MR2172729} and typically rely on symmetric kernels; see, e.g., \cite[p.~43]{MR848134}. However, the traditional location-scale kernel method applies only to data supported on $\R^d$. For other supports, symmetric kernels suffer from boundary bias caused by kernel mass spilling outside the support, which reduces the effective weight of observations near the boundary.

Many strategies have been developed to mitigate this problem. An early approach for handling this issue on the interval $[0, 1]$ is due to \citet{MR564251}, whose estimator is a weighted sum of the response variables whose weights are integrals of a fixed kernel over sets partitioning the input space, each containing a fixed design point. To reduce the bias caused by the fixed kernel, they introduced a \emph{boundary kernel} that solves a specific variational problem near the boundary. The Gasser--M\"uller (GM) estimator is known to outperform the classical Nadaraya–Watson (NW) estimator \citep{MR166874, MR185765} and the Priestley--Chao estimator \citep{MR331616}. This approach was further explored in \cite{MR816088, MR1130920} for the univariate case, in \cite{MR960887, MR1369600, MR1239102} for the multivariate case, and expanded upon in \cite{doi:10.1007/BF00147776, MR1649872, MR1752313} for density estimation.

Among the early alternatives to the GM estimator, the class of local polynomial estimators is particularly prominent; these are constructed by solving locally weighted least-squares regression problems. This method was introduced in \cite{MR443204} and further studied, e.g., in \cite{MR556476, MR582402, MR594650, MR673642}; refer to \cite[Section~3.8]{MR1383587} for bibliographic notes. This approach has the clear advantage of naturally avoiding boundary effects \citep{MR1209561, MR1212173}, without requiring boundary-specific kernel modifications, which simplifies implementation. In 1992, the local linear (LL) estimation method was further refined when \citet{MR1193323} introduced the idea of a variable smoothing parameter that is optimal under the mean integrated squared error (MISE) criterion, allowing the kernel to adapt locally at each point of the support; see \cite{MR1463570} for subsequent refinements. This idea was extended to local polynomial estimators and the multivariate setting in \cite{MR1311979}; see also \cite[Chapter~3]{MR1383587} for a book treatment.

Another alternative is provided by Bernstein regression estimators, studied in the fixed-design and random-design settings in \cite{MR858109} and \cite{MR1437794}, respectively. A general strategy relies on the classical Weierstra{\ss} approximation theorem, which asserts that every continuous function on a closed interval admits a uniform approximation by polynomials. Among the most prominent constructive realizations of this theorem are Bernstein polynomials \cite{Bernstein_1912}, which yield explicit approximations and play a central role in approximation theory. The application of Bernstein polynomials to statistics dates back at least to \cite{MR397977}; see also \cite{MR1910059, MR2270097, MR3474765, MR1685301, MR638651, MR2662607, MR2960952, MR4287788}.

Around the turn of the 21st century, several asymmetric kernels were proposed as alternatives to classical symmetric kernels. The most commonly used asymmetric kernels are the beta kernel on $[0,1]$ and the gamma kernel on $[0, \infty)$. When dealing with $d$-variate data, one can consider products of univariate kernels on product spaces \citep{MR2568128, MR3760293} or kernels that are tailored to non-product spaces such as Dirichlet kernels on the simplex \citep{MR4544604, arXiv:2506.08816, MR4319409}, Wishart kernels on the cone of positive definite matrices \citep{arXiv:2512.08232, MR4358612}, inverse Gaussian kernels on half-spaces \citep{MR4939549}, etc. In the specific case of univariate data, asymmetric kernel-based NW, LL, and GW estimators were discussed, e.g., in \cite{MR1742101, MR1910175, MR4302589, MR4899357, MR4471289, MR4556820, MR3494026, MR3499725}.

Only recently have Dirichlet-kernel estimators been considered for regression surfaces on the simplex. A Dirichlet-kernel version of the NW estimator (D-NW) was introduced and studied in \cite{MR4796622} within the broader framework of conditional $U$-statistics. Similarly, a Dirichlet-kernel version of Fan's LL estimator (D-LL) with a variable smoothing parameter was proposed and studied in \cite{MR4905615}, thereby extending some of the one-dimensional results in \cite{MR1910175}.

This paper has two objectives. First, a Dirichlet-kernel Gasser--M\"uller (D-GM) estimator for fixed design regression on the simplex is introduced and its asymptotic properties are studied, generalizing the work of \citet{MR1742101} in the one-dimensional case. Second, its small-sample performance is compared with that of the D-NW/LL estimators \citep{MR4796622, MR4905615} through simulations. This comparison is motivated by the need to clarify the relative merits of Dirichlet-kernel regression procedures on the simplex. Indeed, the D-GM estimator is a natural analog of the classical Gasser--M\"uller construction in a setting where Dirichlet-kernel methods have only recently been developed, and its finite-sample behavior is not apparent a priori from asymptotic considerations. In particular, for $d\in\{1,2,3\}$, the leading terms of the mean squared error (MSE) and MISE coincide with those of the D-NW/LL estimators, so it is of independent interest to determine whether this asymptotic similarity persists in practice.

The rest of the paper is organized as follows. Section~\ref{sec:definitions.notations} contains preliminary definitions and notations. The assumptions required for the derivations are then presented and discussed in Section~\ref{sec:assumptions}. The main results are stated in Section~\ref{sec:main.results}. The results of a simulation study are then reported in Section~\ref{sec:simulations}, where it is seen that the D-LL estimator performs best across the board. As an illustration, this estimator is applied in Section~\ref{sec:application} to the GEMAS dataset \citep{doi:10.1016/j.scitotenv.2012.02.032}, concerned with the chemical composition of soil samples from various agricultural and grazing lands in Europe. The proofs of the main results are relegated to \ref{app:proofs}, and some technical lemmas on which they depend can be found in \ref{app:tech.lemmas}. A link to the \textsf{R} code that generated the figures, the simulation results and the real-data application is given at the end.

\newpage
\section{Definitions and notations}\label{sec:definitions.notations}

In what follows, the dimension $d \in \N = \{1, 2, \dots\}$ is fixed. The $d$-dimensional simplex and its interior are defined by
\[
\mathcal{S}_d = \big\{\bb{s} \in [0,1]^d: \|\bb{s}\|_1 \leq 1\big\}, \qquad \mathrm{Int}(\mathcal{S}_d) = \big\{\bb{s}\in (0,1)^d: \|\bb{s}\|_1 < 1\big\},
\]
where $\|\bb{s}\|_1 = |s_1| + \cdots + |s_d|$ is the $\ell^1$ norm on $\R^d$. For $\bb{\alpha} = (\alpha_1,\dots,\alpha_d) \in (0, \infty)^d$ and $\beta\in (0,\infty)$, the density of the $\mathrm{Dirichlet}\hspace{0.2mm}(\bb{\alpha},\beta)$ distribution is defined by
\[
K_{\bb{\alpha},\beta}(\bb{s})
= \frac{\Gamma(\|\bb{\alpha}\|_1 + \beta)}{\Gamma(\beta) \prod_{j=1}^d \Gamma(\alpha_j)} \, (1 - \|\bb{s}\|_1)^{\beta - 1} \, \prod_{j=1}^d s_j^{\alpha_j - 1},
\qquad \bb{s} = (s_1,\dots,s_d)\in \mathcal{S}_d.
\]

Let $Y_1, \dots, Y_n$ be the response variables associated with known and fixed design points $\bb{x}_1, \dots, \bb{x}_n\in \mathrm{Int}(\mathcal{S}_d)$. The design density $f$ corresponds to the density of the design points in the limit as $n \to \infty$. More formally, let $\smash{\bb{x}_1^{(n)}, \dots, \bb{x}_n^{(n)}}\in \mathrm{Int}(\mathcal{S}_d)$ denote the expanding grid of design points at sample~size~$n$. It is assumed that, for any Borel set $B\in \mathcal{B}(\R^d)$,
\[
\lim_{n \to \infty} \frac{1}{n} \sum_{i=1}^n \ind_B(\bb{x}_i^{(n)}) = \int_B f(\bb{x}) \, \rd \bb{x}.
\]
Throughout the paper, the simplified notation $\bb{x}_1, \dots, \bb{x}_n$ is used instead of $\bb{x}_1^{(n)}, \dots, \bb{x}_n^{(n)}$ for the design points at sample size $n$.

Assume that, for every $i \in [n] = \{1, \dots, n\}$, $Y_i$ follows the heteroscedastic model
\begin{equation}\label{eq:model}
Y_i = m(\bb{x}_i) + \e_i,
\end{equation}
where $m: \mathcal{S}_d\to \R$ is an unknown regression function and the random error terms $\e_1, \dots, \e_n$ are mutually independent and such that, for each $i \in [n]$, $\e_i$ has mean zero and variance $\sigma^2(\bb{x}_i)$.

Given a smoothing parameter $b \in (0, \infty)$ and a sequence $B_1, \dots, B_n$ of convex compact sets that partitions the simplex $\mathcal{S}_d$ and satisfies $\bb{x}_i\in \mathrm{Int}(B_i)$ for every $i \in [n]$, the D-GM estimator is defined by
\begin{equation}
\label{eq:GM.estimator}
\hat{m}_{n,b}^{\mathrm{GM}}(\bb{s}) = \sum_{i=1}^n Y_i \int_{B_i} K_{\bb{s} / b + \bb{1}, (1 - \|\bb{s}\|_1) / b + 1}(\bb{x}) \, \rd \bb{x}, \qquad \bb{s}\in \mathcal{S}_d,
\end{equation}
where $\bb{1} = (1, \dots, 1)^\top$ is a $d \times 1$ vector of ones. For the sake of simplicity, set $\kappa_{\bb{s},b} = K_{\bb{s} / b + \bb{1}, (1 - \|\bb{s}\|_1) / b + 1}$ for all $\bb{s} \in \mathcal{S}_d$ and $b \in (0, \infty)$, so that the D-GM estimator can be expressed equivalently as
\[
\hat{m}_{n,b}^{\mathrm{GM}}(\bb{s})  = \sum_{i=1}^n Y_i \int_{B_i} \kappa_{\bb{s},b} (\bb{x}) \, \rd \bb{x}.
\]

This estimator is a multivariate analog of the type-3 estimator originally defined by \citet[p.~27]{MR564251}, where the fixed weight function $w$ is replaced by a Dirichlet kernel whose parameters locally adapt with the position of the estimation point~$\bb{s}\in \mathcal{S}_d$. In particular, the Dirichlet-kernel parameters $\bb{\alpha} = \bb{s} / b + \bb{1}$ and $\beta = (1 - \|\bb{s}\|_1) / b + 1$ are chosen so that, as $b \to 0$, the mean and variance of the Dirichlet kernel are asymptotically equal to $\bb{s}$ and zero, respectively. As a result, the kernel concentrates more and more around $\bb{s}$ as the smoothing parameter shrinks. These properties were established in \cite[p.~7]{MR4319409}. Note that the mode is exactly equal to $\bb{s}$.

Throughout the paper, the following conventions are adopted. For any given set $E$, $\mathrm{Int}(E)$ denotes the interior of $E$ and $\partial E$ denotes the boundary of $E$. The notation $u = \OO(v)$ means that $\limsup |u / v| \leq C < \infty$ as $n \to \infty$ or $b\to 0$, depending on the context. The positive constant $C$ may depend on the regression function $m$, the variance function $\sigma^2$, the design density $f$, and the dimension $d$, but on no other variables unless explicitly written as a subscript. A common occurrence is a local dependence of the asymptotics on a given estimation point~$\bb{s}\in \mathcal{S}_d$, in which case one writes $u = \OO_{\bb{s}}(v)$. The alternative notation $u \ll v$ is also used to mean $u, v\in [0, \infty)$ and $u = \OO(v)$. If both $u \ll v$ and $u \gg v$ hold, one writes $u \asymp v$. Similarly, the notation $u = \oo(v)$ means that $\lim |u / v| = 0$ as $n \to \infty$ or $b\to 0$. Subscripts indicate which parameters the convergence rate can depend on. The symbol $\rightsquigarrow$ denotes convergence in distribution, $\lambda$ denotes the Lebesgue measure in $\R^d$, and the shorthand notations $[d] = \{1, \dots, d\}$ and $[n] = \{1, \dots, n\}$ are used frequently. Finally, note that the smoothing parameter $b = b(n)$ is always implicitly a function of $n$, the number of observations, except for Lemmas~\ref{lem:local.bound}~and~\ref{lem:A.b.asymptotics} in \ref{app:tech.lemmas}.

\section{Assumptions}\label{sec:assumptions}

The assumptions used to establish the results stated in Section~\ref{sec:main.results} are the following.

\vspace{-1mm}
\begin{enumerate}[label=(A\arabic*)]\setlength\itemsep{0em}
\item $f$ and $\sigma^2$ are Lipschitz continuous on $\mathcal{S}_d$, and $m$ is twice continuously differentiable on $\mathcal{S}_d$. In particular, $\sigma^2$ is bounded on $\mathcal{S}_d$, i.e., $\max\{\sigma^2(\bb{x}) : \bb{x}\in \mathcal{S}_d\} < \infty$. \label{ass:1}
\item $\min \{ f(\bb{x}) : \bb{x}\in \mathcal{S}_d \} > 0$ and $\min \{ \sigma^2(\bb{x}) : \bb{x}\in \mathcal{S}_d \} > 0$. \label{ass:2}
\item $b = b(n)\to 0$ and $n^{-1/d} = \oo(b)$ as $n\to \infty$. \label{ass:3}
\item For every fixed $n \in \N$, there exist convex compact subsets $B_1, \dots, B_n$ of $\mathcal{S}_d$ such that, \label{ass:4}
\begin{enumerate}[label=(A4.\alph*)]\setlength\itemsep{0em}
\item $\mathrm{Int}(B_i) \cap \mathrm{Int}(B_j) = \emptyset$ for all $i\neq j$, and $B_1 \cup \cdots \cup B_n = \mathcal{S}_d$; \label{ass:4.a}
\item $\bb{x}_i\in \mathrm{Int}(B_i)$ for every $i\in [n]$; \label{ass:4.b}
\item $\lambda(\partial B_i) = 0$, for every $i\in [n]$; \label{ass:4.c}
\item $\mathrm{diam}(B_i) = \OO(n^{-1/d})$ uniformly in $i$, as $n\to\infty$, where $\mathrm{diam}(B_i)$ is the diameter of~$B_i$; \label{ass:4.d}
\item $|\int_{B_i} f(\bb{t}) \, \rd \bb{t}-n^{-1} | = \OO(n^{-1-1/d})$ uniformly in $i$, as $n \to \infty$. \label{ass:4.e}
\end{enumerate}
\end{enumerate}

Each of these assumptions is discussed in turn, starting with Assumption~\ref{ass:4}, which is the most crucial to the D-GM estimator.

Assumption~\ref{ass:4.a} implies that, for any fixed $n \in \N$, the finite sequence $B_1, \dots, B_n$ of convex compact sets partitions the simplex $\mathcal{S}_d$. For any given estimation point $\bb{s}\in \mathcal{S}_d$, the D-GM estimator is a linear combination of the response variables, whose weights are controlled by the Dirichlet kernel $\kappa_{\bb{s},b}$, which is designed to adapt its shape on the simplex depending on the position of $\bb{s}$, thereby providing asymptotically negligible bias near the boundary.

Assumption~\ref{ass:4.b} means that, for every integer $i \in [n]$, the design point $\bb{x}_i$ is a representative of the set $B_i$. Their role is analogous to the points at which the height of a function is evaluated within the sets of a partition in a Riemann sum.

Assumption~\ref{ass:4.c}, together with Assumption~\ref{ass:4.a}, ensures that the integral of any integrable function over $\mathcal{S}_d$ can be decomposed into the sum of integrals over the respective sets $B_1, \dots, B_n$. This decomposition of the integral is useful, e.g., for the pointwise bias analysis in the proof of Proposition~\ref{prop:bias}, and for the pointwise variance analysis in the proof of Proposition~\ref{prop:var}; see Appendix~\ref{Appendix:A1} and Appendix~\ref{Appendix:A2}, respectively.

Assumption~\ref{ass:4.d}, which is technical, provides an asymptotic upper bound on the diameter of the sets $B_1, \dots, B_n$. This assumption is used repeatedly in the proofs, together with the smoothness conditions in Assumption~\ref{ass:1}, to control differences of $f$, $\sigma^2$, or $m$, at two different points within each set. Such differences appear in multiple places in the pointwise variance analysis when applying the multivariate mean value theorem.

\begin{figure}[t!]
\centering
\includegraphics[width=0.40\textwidth]{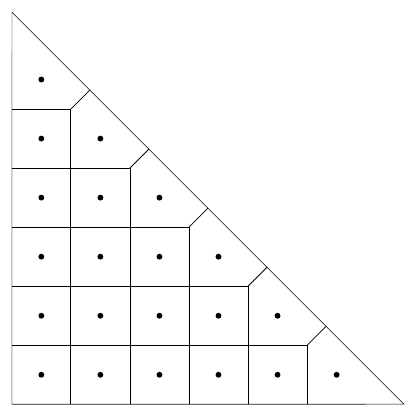}
\caption{The black dots represent the design points $\bb{x}_1, \dots, \bb{x}_n$ and are chosen here to form a grid of mesh size $\asymp n^{-1/2}$ over the two-dimensional simplex $\mathcal{S}_2$. For each $i\in [n]$, the polygonal set $B_i$ is the Voronoi cell of $\bb{x}_i$.}
\label{fig:Voronoi}
\end{figure}

Assumption~\ref{ass:4.e} is another technical requirement which imposes a uniform control on the relative weight of the sets $B_1, \dots, B_n$. For example, if the design density is uniform, then Assumption~\ref{ass:4.e} indicates that each set $B_i$ has asymptotic weight $1/n$, which is natural. This would be the case, for instance, when the design points form a grid of mesh size $\asymp n^{-1/d}$ over $\mathcal{S}_d$ and the sets $B_1, \dots, B_n$ are chosen to be the associated Voronoi cells; see Figure~\ref{fig:Voronoi} for the case $d = 2$.

Note that if random design points are uniformly generated for every integer $n\in \N$, with each sequence fixed once it has been generated, and if for each integer $i\in [n]$, the set $B_i$ is taken to be the Voronoi cell of $\bb{x}_i$, then the asymptotic upper bound on the diameter of $B_i$ in Assumption~\ref{ass:4.d} would need to be relaxed by a logarithmic factor, as shown in \cite[Theorem~5.1]{MR3668473}; see also \cite{MR4123649} for a closely related analysis. In turn, this additional factor in Assumption~\ref{ass:4.d} would have an impact on the asymptotics of every result in Section~\ref{sec:main.results}.

Assumption~\ref{ass:1} is common in nonparametric regression contexts; see, e.g., \cite[Chapter~3]{MR1383587}. The requirement that the regression function $m$ is twice continuously differentiable is needed to apply a second-order Taylor expansion in the pointwise bias analysis in the proof of Proposition~\ref{prop:bias}. The assumption that both $f$ and $\sigma^2$ are Lipschitz continuous means that for each integer $i\in [n]$, the value of $f$ or $\sigma^2$ at any given point within each set $B_i$ is equal to its value at the representative $\bb{x}_i\in B_i$, up to an error asymptotically bounded by the distance from the point to $\bb{x}_i$. More generally, Assumption~\ref{ass:1} implies that for each integer $i\in [n]$, differences of $f$ or $\sigma^2$ at two different points in $B_i$ are asymptotically bounded by the diameter of $B_i$. Hence, this assumption is often employed in the proofs in conjunction with Assumption~\ref{ass:4.d}.

Assumption~\ref{ass:2} is also standard; see, e.g., \cite[Chapter~3]{MR1383587}. This assumption is used to control the ratio $\sigma^2/ f$ and related error terms, which appear in the asymptotics of the variance of the D-GM estimator; see the proof of Proposition~\ref{prop:var}. Additionally, it is indispensable to apply Lebesgue's dominated convergence theorem (i.e., exchange limits and integrals) when investigating the integrated variance in the proof of Theorem~\ref{thm:MISE}.

Finally, Assumption~\ref{ass:3} states not only that $b = b(n)$ is a function of the sample size $n$, as emphasized at the end of Section~\ref{sec:definitions.notations}, but that both $b = \oo(1)$ and $n^{-1/d} = \oo(b)$ as $n \to \infty$. The latter condition, which is crucial, will be used repeatedly to handle various error terms in the proofs. It implies that for each integer $i\in [n]$, the diameter of the set $B_i$ (for example $\asymp n^{-1/d}$ with $d = 2$ in Figure~\ref{fig:Voronoi}) must be asymptotically smaller than the smoothing parameter $b$ as $n \to \infty$. In particular, note that if the smoothing parameter $b$ is taken to be proportional to $n^{-2/(d+4)}$, the dimension restriction $d\in \{ 1, 2, 3 \}$ must then be imposed to obtain certain theoretical guarantees on the asymptotic rates of the MSE and MISE (see Corollary~\ref{cor:MSE} and Theorem~\ref{thm:MISE}), which is an obvious limitation of the D-GM estimator. See Remark~\ref{rem:optim.MSE} for additional details.

\section{Main results}\label{sec:main.results}

The pointwise bias/variance of the one-dimensional beta-GM estimator on $[0,1]$ were first studied in \cite{MR1742101}. Propositions~\ref{prop:bias} and~\ref{prop:var} below generalize these findings to the D-GM estimator when the domain is the simplex, for an arbitrary $d\in \N$.

\begin{proposition}[Pointwise bias]\label{prop:bias}
Suppose that Assumptions~\ref{ass:1}--\ref{ass:4} hold. Then, as $n \to \infty$ and uniformly in $\bb{s}\in \mathcal{S}_d$, one has
\[
\Bias\{\hat{m}_{n,b}^{\mathrm{GM}}(\bb{s})\} = \EE\{\hat{m}_{n,b}^{\mathrm{GM}}(\bb{s})\} - m(\bb{s}) = b \, g(\bb{s}) + \oo(b),
\]
where
\begin{equation}\label{eq:fonction-g}
g(\bb{s}) = \sum_{j=1}^d \{1 - (d+1) s_j \} \, \frac{\partial}{\partial s_j} \, m(\bb{s}) + \frac{1}{2} \sum_{j,k=1}^d s_j (\ind_{\{j = k\}} - s_k) \, \frac{\partial^2}{\partial s_j \partial s_k} \, m(\bb{s}).
\end{equation}
\end{proposition}

\begin{remark}
The leading term in the pointwise bias of the D-GM estimator is the same as that of the D-NW/LL estimators \cite{MR4796622,MR4905615}. As in \cite[Remark~1]{MR4905615}, it is conjectured that the D-LL estimator exhibits improved boundary behavior due to the influence of higher-order terms in its asymptotic expansion, which may help explain its superior performance in the simulation study presented in Section~\ref{sec:simulations}.
\end{remark}

\begin{proposition}[Pointwise variance]\label{prop:var}
Suppose that Assumptions~\ref{ass:1}--\ref{ass:4} hold, and let $\mathcal{J} \subseteq [d]$ and $\bb{\lambda} = (\lambda_1,\dots,\lambda_d)\in (2,\infty)^d$ be given. If a sequence $\bb{s} = \bb{s}(b)$ is selected inside $\mathrm{Int}(\mathcal{S}_d)$ such that $s_j / b \to \lambda_j$ for every integer $j \in \mathcal{J}$, and $s_j$ is fixed for every integer $j \in [d] \setminus \mathcal{J}$, then, as $n \to \infty$,
\[
\Var\{\hat{m}^{\mathrm{GM}}_{n,b}(\bb{s})\}
= n^{-1} b^{-(d + |\mathcal{J}|)/2} \, \left\{\frac{\psi_{\mathcal{J}}(\bb{s}) \sigma^2(\bb{s})}{f(\bb{s})} \prod_{j\in \mathcal{J}} \frac{\Gamma(2\lambda_j + 1)}{2^{2\lambda_j + 1} \Gamma^2(\lambda_j + 1)} + \oo_{\bb{\lambda},\bb{s}}(1)\right\},
\]
where, for every subset $\mathcal{J} \subseteq [d]$ and vector $\bb{s}\in \mathrm{Int}(\mathcal{S}_d)$,
\[
\psi_{\mathcal{J}}(\bb{s}) = \left\{(4\pi)^{d - |\mathcal{J}|} \, (1 - \|\bb{s}\|_1) \prod_{j \in [d] \backslash\mathcal{J}} s_j \right\}^{-1/2}.
\]
\end{proposition}

\begin{remark}
In view of Proposition~\ref{prop:var}, the pointwise variance of the D-GM estimator at any fixed estimation point $\bb{s}\in \mathrm{Int}(\mathcal{S}_d)$ has an asymptotic order of $n^{-1} b^{-d/2}$. When $\bb{s}$ approaches the boundary in one or more coordinates, the variance increases by a multiplicative factor of $b^{-1/2}$ per coordinate. Thus, if $\bb{s}$ is near a $(d - |\mathcal{J}|)$-dimensional face of the boundary, the pointwise variance scales as $n^{-1} b^{-(d + |\mathcal{J}|)/2}$. The D-NW/LL estimators \cite{MR4796622,MR4905615} share the same features; see \cite[Proposition~1]{MR4905615}.
\end{remark}

As an immediate corollary to Propositions~\ref{prop:bias} and \ref{prop:var}, the MSE is provided below for any fixed estimation point $\bb{s}\in \mathrm{Int}(\mathcal{S}_d)$.

\begin{corollary}[MSE; interior point]\label{cor:MSE}
Suppose that Assumptions~\ref{ass:1}--\ref{ass:4} hold, and let $\bb{s}\in \mathrm{Int}(\mathcal{S}_d)$ be fixed. Then, as $n \to \infty$,
\[
\begin{aligned}
\mathrm{MSE}\{\hat{m}_{n,b}^{\mathrm{GM}}(\bb{s})\}
&= \big[\Bias\{\hat{m}_{n,b}^{\mathrm{GM}}(\bb{s})\}\big]^2 + \Var\{\hat{m}_{n,b}^{\mathrm{GM}}(\bb{s})\} \\
&= b^2 g^2(\bb{s}) + n^{-1} b^{-d/2} \, v(\bb{s}) + \oo(b^2) + \oo_{\bb{s}}(n^{-1} b^{-d/2}),
\end{aligned}
\]
where $v(\bb{s}) = \psi_{\emptyset}(\bb{s}) \sigma^2(\bb{s}) / f(\bb{s})$. In particular, if $g(\bb{s}) \sigma^2(\bb{s}) \neq 0$ and $d\in \{1,2,3\}$, the asymptotically (pointwise) optimal smoothing parameter $b$, which minimizes $b^2g^2(\bb{s}) + n^{-1}b^{-d/2}v(\bb{s})$, is given by
\[
b_{\mathrm{opt}}(\bb{s}) = n^{-2/(d+4)} \left\{\frac{d}{4} \, \frac{v(\bb{s}) }{ g^2(\bb{s})}\right\}^{2/(d+4)},
\]
where
\[
\mathrm{MSE}\{\hat{m}_{n,b_{\mathrm{opt}}(\bb{s})}^{\mathrm{GM}}(\bb{s})\}
= n^{-4 / (d+4)} \left\{\frac{1 + d/4}{(d/4)^{d/(d+4)}}\right\} \big\{v(\bb{s})\big\}^{4 / (d+4)} \big\{g^2(\bb{s})\big\}^{d/(d+4)} + \oo_{\bb{s}} \big(n^{-4/(d+4)} \big).
\]
\end{corollary}

\begin{remark}\label{rem:optim.MSE}
The restriction $d\in \{1,2,3\}$ ensures that choosing $b_{\mathrm{opt}}(\bb{s}) \asymp n^{-2 / (d+4)}$ is feasible in view of Assumption~\ref{ass:3}. Indeed, $\smash{n^{-1/d} = \oo\{b_{\mathrm{opt}}(\bb{s})\}}$ is equivalent to $n^{2 / (d+4)} = \oo(n^{1/d})$, i.e., $2 / (d+4) < 1 / d$, or equivalently $d < 4$. The same considerations apply to Theorem~\ref{thm:MISE} below.
\end{remark}

The following result provides the MISE of the D-GM estimator, showing that the contribution from points near the boundary is negligible. A similar result was derived for the classical multivariate GM estimators with $k$th-order kernels in \cite[p.~179--180]{MR1369600}; see also \cite[Chapter~6]{MR960887} and~\cite{MR1239102}.

\begin{theorem}[MISE]\label{thm:MISE}
Suppose that Assumptions~\ref{ass:1}--\ref{ass:4} hold. Then, as $n \to \infty$,
\begin{align*}
\mathrm{MISE}(\hat{m}_{n,b}^{\mathrm{GM}}) & = \int_{\mathcal{S}_d} \mathrm{MSE}\{\hat{m}_{n,b}^{\mathrm{GM}}(\bb{s})\} \, \rd \bb{s} \\
& = b^2 \int_{\mathcal{S}_d} g^2(\bb{s}) \, \rd \bb{s} + n^{-1} b^{-d/2} \int_{\mathcal{S}_d} v(\bb{s}) \, \rd \bb{s} + \oo(b^2 + n^{-1} b^{-d/2}).
\end{align*}
In particular, if $g \not\equiv 0$, $\sigma^2 \not\equiv 0$ and $d\in \{1,2,3\}$, the asymptotically (global) optimal smoothing parameter $b$, which minimizes $b^2 \int_{\mathcal{S}_d} g^2(\bb{s}) \, \rd \bb{s} + n^{-1} b^{-d/2} \int_{\mathcal{S}_d} v(\bb{s}) \, \rd \bb{s}$, is given by
\[
b_{\mathrm{opt}} = n^{-2/(d+4)} \left\{ \frac{d}{4} \, \frac{\int_{\mathcal{S}_d} v(\bb{s}) \, \rd \bb{s}}{\int_{\mathcal{S}_d} g^2(\bb{s}) \, \rd \bb{s}}\right\}^{2/(d+4)},
\]
where
\[
\mathrm{MISE}(\hat{m}_{n,b_{\mathrm{opt}}}^{\mathrm{GM}})
= n^{-4 / (d+4)} \left\{\frac{1 + d/4}{(d/4)^{d/(d+4)}}\right\} \frac{\big\{\int_{\mathcal{S}_d}v(\bb{s}) \, \rd \bb{s}\big\}^{4 / (d+4)}}{\big\{\int_{\mathcal{S}_d} g^2(\bb{s}) \, \rd \bb{s}\big\}^{-d/(d+4)}} + \oo(n^{-4/(d+4)}).
\]
\end{theorem}

\begin{remark}
Propositions~\ref{prop:bias} and \ref{prop:var}, Corollary~\ref{cor:MSE}, and Theorem~\ref{thm:MISE} remain valid under the weaker assumption that the error terms $\e_1, \dots, \e_n$ are uncorrelated. Indeed, the independence is only used to justify the first equality of Eq.~\eqref{eq:variance.expansion} in the proof of Proposition~\ref{prop:var}.
\end{remark}

Verification of the Lindeberg condition for triangular arrays leads to the asymptotic normality of the D-GM estimator at any fixed estimation point $\bb{s}\in \mathrm{Int}(\mathcal{S}_d)$. For the asymptotic normality of the classical multivariate GM estimators with $k$th-order kernels, refer to \cite{MR1239102}.

\begin{theorem}[Asymptotic normality]\label{thm:CLT}
Let $d\in \{1,2,3\}$. Suppose that Assumptions~\ref{ass:1}--\ref{ass:4} hold, and consider a fixed point $\bb{s}\in \mathrm{Int}(\mathcal{S}_d)$ such that $\sigma^2(\bb{s})\in (0,\infty)$. Let $R_n = n^{1/2} b^{d/4}$, and assume that the errors $\e_1, \dots, \e_n$ satisfy
\begin{equation}\label{eq:asympotic.normality.condition}
\forall_{\delta\in (0,\infty)} \qquad \lim_{n\to \infty} \frac{1}{R_n^2} \sum_{i=1}^n \EE\big(\e_i^2 \, \ind_{\{|\e_i| > \delta R_n\}}\big) = 0.
\end{equation}
Then the following statements hold true.
\begin{itemize}\setlength\itemsep{0em}
\item [(i)]
If $R_n b \to 0$ as $n \to \infty$, then $R_n \{\hat{m}_{n,b}^{\mathrm{GM}}(\bb{s}) - m(\bb{s})\} \rightsquigarrow \mathcal{N} \big[ 0, v(\bb{s}) \big]$.
\item [(ii)]
If $n^{2/(d+4)} \, b \to \nu$ as $n \to \infty$ for some $\nu\in (0,\infty)$, then
\[
n^{2 / (d+4)} \{\hat{m}_{n,b}^{\mathrm{GM}}(\bb{s}) - m(\bb{s})\} \rightsquigarrow \mathcal{N} \big[ \nu \, g(\bb{s}), \nu^{-d/2} v(\bb{s}) \big].
\]
\end{itemize}
\end{theorem}

\begin{remark}
Condition~\eqref{eq:asympotic.normality.condition} is implied by the Lyapunov condition:
\[
\exists_{a\in (0,\infty)} \qquad \lim_{n\to \infty} \frac{1}{R_n^{2+\delta}} \sum_{i=1}^n \EE(|\e_i|^{2+a}) = 0;
\]
see, e.g., \citet[p.~362]{MR1324786}.
\end{remark}

\section{Simulation study}\label{sec:simulations}

This section describes the results of a modest simulation study designed to investigate the small-sample performance of the proposed D-GM estimator. Its competitors are the D-NW estimator \cite{MR4796622} and the D-LL estimator \cite{MR4905615}, respectively given by
\[
\hat{m}^{\mathrm{NW}}_{n,b}(\bb{s}) = \frac{\sum_{i=1}^n Y_i \, \kappa_{\bb{s},b}(\bb{x}_i)}{\sum_{i=1}^n \kappa_{\bb{s},b}(\bb{x}_i)}
\qquad \mbox{and} \qquad
\hat{m}^{\mathrm{LL}}_{n,b}(\bb{s}) = \hat{\alpha}_{\bb{s}} = \bb{e}_1^{\top} (\mathcal{X}_{\bb{s}}^{\top} W_{\bb{s}} \mathcal{X}_{\bb{s}})^{-1} \mathcal{X}_{\bb{s}}^{\top} W_{\bb{s}} \bb{Y},
\]
where $\bb{e}_1 = (1, 0, \dots, 0)^\top$ is a $(d+1) \times 1$ vector, and
\[
\mathcal{X}_{\bb{s}} =
\begin{bmatrix}
1 & (\bb{x}_1 - \bb{s})^{\top} \\
\vdots & \vdots \\
1 & (\bb{x}_n - \bb{s})^{\top}
\end{bmatrix}_{n \times (d+1)}, \qquad
W_{\bb{s}} = \mathrm{diag}
\begin{bmatrix}
\kappa_{\bb{s},b}(\bb{x}_1) \\
\vdots \\
\kappa_{\bb{s},b}(\bb{x}_n)
\end{bmatrix}_{n \times 1}, \qquad
\bb{Y} =
\begin{bmatrix}
Y_1 \\
\vdots \\
Y_n
\end{bmatrix}_{n \times 1}.
\]
The D-LL estimator corresponds to the intercept $\hat{\alpha}_{\bb{s}}$ that minimizes the locally weighted loss
\[
L(\alpha,\bb{\beta}) = \sum_{i=1}^n \{Y_i - \alpha - \bb{\beta}^{\top} (\bb{x}_i - \bb{s})\}^2 \kappa_{\bb{s},b}(\bb{x}_i).
\]
The design points $\bb{x}_1, \dots, \bb{x}_n\in \smash{\mathrm{Int}(\mathcal{S}_d)}$ were chosen according to the fixed mesh in Figure~\ref{fig:Voronoi}, viz.
\[
\{(w_k (i - 1) + 1/2, w_k (k - j) + 1/2)/(k+1) : 1 \leq i \leq j \leq k \},
\]
where $w_k = (k - 1 / \sqrt{2}) / (k - 1)$ and $k$ is a positive integer. In the case of the D-GM estimator, the sequence $B_1, \dots, B_n$ is chosen to be the Voronoi diagram associated with the design points, so that, for every $i\in [n]$, $B_i$ is the Voronoi cell of $\bb{x}_i$. For each integer $i\in [n]$, the response variable $Y_i$ was generated from Model~\eqref{eq:model} with
\[
\e_1, \dots, \e_n\stackrel{\mathrm{iid}}{\sim} \mathcal{N} \left[ 0, \frac{1}{10} \, \mathrm{IQR}\{m(\bb{x}_1), \dots, m(\bb{x}_n)\} \right],
\]
where $\mathrm{IQR}$ denotes the interquartile range. The smoothing parameter was selected by leave-one-out cross-validation (LOOCV). For a given method in $\{\mathrm{D\text{-}GM}, \mathrm{D\text{-}NW}, \mathrm{D\text{-}LL}\}$ and an observed sample $(\bb{x}_i,y_i)_{i=1}^n$, the least-squares criterion is
\[
\mathrm{LOOCV}^{\mathrm{method}}(b) = \frac{1}{n}\sum_{i=1}^n \big\{y_i-\hat{m}_{n,b,(-i)}^{\mathrm{method}}(\bb{x}_i)\big\}^2,
\]
where $\smash{\hat{m}_{n,b,(-i)}^{\mathrm{method}}}$ denotes the corresponding estimator computed without the $i$th observation. Six target regression functions defined for all $\bb{s} = (s_1, s_2)\in \mathcal{S}_2$ were then tested, namely
\[
\begin{array}{lll}
m_1(\bb{s}) = \ln (1 + s_1 + s_2), &
m_2(\bb{s}) = \sin(s_1) + \cos(s_2), &
m_3(\bb{s}) = \sqrt{s_1} + \sqrt{s_2}, \\
m_4(\bb{s}) = s_1 (1 + s_2), &
m_5(\bb{s}) = (s_1 + 1/4)^2 + (s_2 + 3/4)^2, &
m_6(\bb{s}) = (1 + s_1) e^{s_2}.
\end{array}
\]

For each of the three methods, each target regression function, and each sample size $n = k (k+1)/2$ with $k\in \{7, 10, 14\}$, Table~\ref{tab:1} reports the mean, median, standard deviation (SD), and IQR of the sequence of Monte Carlo estimates of the integrated squared errors (ISEs), i.e.,
\[
\widetilde{\mathrm{ISE}}_{\mathrm{method}}^{j,r} = \frac{1}{1000 \times 2} \sum_{i=1}^{1000} \big\{\hat{m}_{n,\hat{b}_{n,r}}^{\mathrm{method}}(\bb{U}_i) - m_j(\bb{U}_i)\big\}^2,
\]
where $r\in \{1, \dots, 100 \}$ indexes the Monte Carlo replications, $\hat{b}_{n,r}$ is the LOOCV-optimized smoothing parameter for the $r$th replication, $\bb{U}_1, \dots, \bb{U}_{1000}$ form a random sample from the uniform distribution on $\mathcal{S}_2$, and the factor $2$ in the denominator is the normalization constant for this distribution. In all cases, the D-LL estimator \cite{MR4905615} has a smaller mean and median than the D-GM/NW estimators. The same is true for SD and IQR except in 2 cases: $(m_2,n=55)$ and $(m_3,n=55)$.

\begin{table}[htbp]
\caption{Comparison of the D-GM, D-NW, and D-LL methods based on the mean, median, standard deviation (SD), and interquartile range (IQR) of 100 $\smash{\widetilde{\mathrm{ISE}}}$ values, multiplied by $10^6$, for regression functions $m_1$ through $m_6$ and sample sizes $n\in \{28, 55, 105\}$. The integrals in the definition of the D-GM estimator are computed using the $\mathsf{R}$ command \texttt{adaptIntegrate} with a relative tolerance of $10^{-3}$.}
\label{tab:1}

\bigskip
\renewcommand{\arraystretch}{0.87} 
\setlength{\tabcolsep}{10pt} 
\centering
{\small
\begin{tabular}{ccccccc}
Target & $n$ & Method & Mean & SD & Median & IQR \\
\toprule
\multirow{3}{*}{} & \multirow{3}{*}{28}
 & D-GM & ~1884 & ~242 & ~1826 & ~332 \\
 & & D-NW & ~1778 & ~192 & ~1753 & ~244 \\
 & & D-LL & ~~397 & ~122 & ~~393 & ~168 \\

\multirow{3}{*}{$m_1$} & \multirow{3}{*}{55}
 & D-GM & ~1306 & ~188 & ~1265 & ~222 \\
 & & D-NW & ~1165 & ~168 & ~1155 & ~232 \\
 & & D-LL & ~~576 & ~119 & ~~568 & ~152 \\

\multirow{3}{*}{} & \multirow{3}{*}{105}
 & D-GM & ~1200 & ~111 & ~1208 & ~150 \\
 & & D-NW & ~~955 & ~111 & ~~961 & ~153 \\
 & & D-LL & ~~644 & ~~89 & ~~644 & ~137 \\
\toprule
\multirow{3}{*}{} & \multirow{3}{*}{28}
 & D-GM & ~4030 & ~831 & ~3943 & 1065 \\
 & & D-NW & ~5659 & ~830 & ~5554 & 1240 \\
 & & D-LL & ~2021 & ~527 & ~1924 & ~766 \\

\multirow{3}{*}{$m_2$} & \multirow{3}{*}{55}
 & D-GM & ~3107 & ~421 & ~3113 & ~530 \\
 & & D-NW & ~3251 & ~477 & ~3232 & ~667 \\
 & & D-LL & ~1723 & ~367 & ~1765 & ~544 \\

\multirow{3}{*}{} & \multirow{3}{*}{105}
 & D-GM & ~3164 & ~297 & ~3178 & ~432 \\
 & & D-NW & ~2559 & ~309 & ~2507 & ~467 \\
 & & D-LL & ~1766 & ~229 & ~1751 & ~323 \\
\toprule
\multirow{3}{*}{} & \multirow{3}{*}{28}
 & D-GM & ~5273 & ~610 & ~5290 & ~763 \\
 & & D-NW & ~4145 & ~536 & ~4129 & ~661 \\
 & & D-LL & ~1614 & ~510 & ~1525 & ~615 \\

\multirow{3}{*}{$m_3$} & \multirow{3}{*}{55}
 & D-GM & ~3281 & ~318 & ~3273 & ~365 \\
 & & D-NW & ~2418 & ~306 & ~2409 & ~418 \\
 & & D-LL & ~1275 & ~314 & ~1208 & ~407 \\

\multirow{3}{*}{} & \multirow{3}{*}{105}
 & D-GM & ~3119 & ~217 & ~3107 & ~245 \\
 & & D-NW & ~1853 & ~243 & ~1819 & ~296 \\
 & & D-LL & ~1304 & ~201 & ~1267 & ~234 \\
\toprule
\multirow{3}{*}{} & \multirow{3}{*}{28}
 & D-GM & ~5644 & 1085 & ~5510 & 1487 \\
 & & D-NW & ~5518 & ~862 & ~5467 & 1358 \\
 & & D-LL & ~2468 & ~722 & ~2482 & 1065 \\

\multirow{3}{*}{$m_4$} & \multirow{3}{*}{55}
 & D-GM & ~3744 & ~521 & ~3727 & ~631 \\
 & & D-NW & ~3458 & ~466 & ~3463 & ~583 \\
 & & D-LL & ~2190 & ~413 & ~2234 & ~577 \\

\multirow{3}{*}{} & \multirow{3}{*}{105}
 & D-GM & ~3281 & ~340 & ~3247 & ~444 \\
 & & D-NW & ~2713 & ~376 & ~2710 & ~545 \\
 & & D-LL & ~2059 & ~293 & ~2020 & ~369 \\
\toprule
\multirow{3}{*}{} & \multirow{3}{*}{28}
 & D-GM & 21132 & 3294 & 21092 & 4523 \\
 & & D-NW & 20072 & 2405 & 19759 & 3230 \\
 & & D-LL & ~6129 & 1906 & ~5964 & 2512 \\

\multirow{3}{*}{$m_5$} & \multirow{3}{*}{55}
 & D-GM & 12915 & 1643 & 12661 & 2311 \\
 & & D-NW & 11594 & 1566 & 11595 & 2070 \\
 & & D-LL & ~5914 & 1289 & ~5905 & 1601 \\

\multirow{3}{*}{} & \multirow{3}{*}{105}
 & D-GM & 10970 & ~908 & 10929 & 1171 \\
 & & D-NW & ~8263 & ~913 & ~8178 & 1350 \\
 & & D-LL & ~5451 & ~758 & ~5403 & 1055 \\
\toprule
\multirow{3}{*}{} & \multirow{3}{*}{28}
 & D-GM & 14012 & 2070 & 13803 & 2881 \\
 & & D-NW & 10338 & 1510 & 10272 & 1866 \\
 & & D-LL & ~3538 & 1098 & ~3540 & 1375 \\

\multirow{3}{*}{$m_6$} & \multirow{3}{*}{55}
 & D-GM & ~8354 & 1052 & ~8271 & 1349 \\
 & & D-NW & ~5847 & ~877 & ~5833 & 1190 \\
 & & D-LL & ~3109 & ~691 & ~3100 & ~826 \\

\multirow{3}{*}{} & \multirow{3}{*}{105}
 & D-GM & ~7034 & ~568 & ~7034 & ~752 \\
 & & D-NW & ~4525 & ~633 & ~4483 & ~856 \\
 & & D-LL & ~3076 & ~488 & ~3043 & ~789 \\
\bottomrule
\end{tabular}}
\end{table}

\begin{remark}\label{rem:speed}
In Table~\ref{tab:1}, the integrals in the definition of the D-GM estimator are computed using the $\mathsf{R}$ command \texttt{adaptIntegrate} with a relative tolerance of $10^{-3}$. Independent simulations (not presented here) using small samples indicate that the performance of the D-GM estimator seems to improve to the level of the D-NW estimator when the relative tolerance of \texttt{adaptIntegrate} is decreased to $10^{-5}$. Unfortunately, the D-GM method requires substantially greater computational time than the other two methods, so it was not feasible to generate an analog of Table~\ref{tab:1} with a relative tolerance of $10^{-5}$.
\end{remark}

\begin{remark}
The implementation of the D-GM estimator using any kernel $k$ requires the numerical integration of $\int_{B_i} k(\bb{x}) \, \rd \bb{x}$ for every $i\in[n]$, unlike the NW/LL estimators. As a consequence, the D-GM estimator is slower, especially in high dimensions, and reducing the relative tolerance of the numerical integrals is even more computationally costly; see Remark~\ref{rem:speed}. Therefore, the D-GM estimator is not recommended in practice, although the asymptotic theory for $d\in \{1,2,3\}$ was elucidated in this paper. Ultimately, these computational drawbacks strongly reinforce our recommendation to utilize the more precise and efficient D-LL estimator for nonparametric regression on the simplex.
\end{remark}

\section{Real-data application}\label{sec:application}

The GEMAS (Geochemical mapping of agricultural and grazing land soils) project is a large-scale geochemical survey conducted across Europe, managed by EuroGeoSurveys \citep{doi:10.1016/j.scitotenv.2012.02.032}. The primary aim of this project was to analyze the chemical composition of soil samples from various agricultural and grazing lands. By establishing a consistent and systematic sampling strategy, the project gathered a comprehensive dataset spanning 33 European countries. This effort allowed for the identification of baseline concentrations of major and trace elements in topsoil, providing insights into natural and anthropogenic influences on soil composition.

Each soil sample was collected from the top 20 cm of the soil profile to standardize the dataset. The GEMAS dataset, accessible in the GitHub repository \cite{DaayebGenestKhardaniKlutchnikoffOuimet2026github}, is notable for its broad coverage, extensive number of samples, and rigorous analytical techniques. Among the suite of data collected are properties such as elemental concentrations, particle size distribution, and pH levels measured in a calcium chloride (CaCl$_2$) solution. The use of CaCl$_2$ provides a stabilized ionic strength, which results in more consistent and realistic pH readings, better representing the conditions encountered by plant roots in soil.

After the few rows with missing information were removed, the dataset comprises 2083 design points of the form $\bb{x}_i = (x_{i,1}, x_{i,2})$, which represent the (renormalized) proportions of sand and silt in each soil sample, respectively. The proportion of clay is determined by the complement, $1 - x_{i,1} - x_{i,2}$. The pH in CaCl$_2$ of each sample, denoted $y_i$, is the response variable. The goal is to use the D-LL estimator to estimate the pH in CaCl$_2$ based on soil composition. This nonparametric regression estimate helps in understanding how soil texture influences pH levels, which is critical for agricultural management and assessing soil quality.

The simulation study reported in Section~\ref{sec:simulations} suggests that the D-LL estimator outperforms both the new D-GM estimator and the D-NW estimator. For this reason, only the D-LL estimator is considered in what follows. Note, however, that the other two estimators would lead essentially to the same conclusion in this specific application.

\begin{figure}[t!]
\centering
\includegraphics[width=0.44\textwidth]{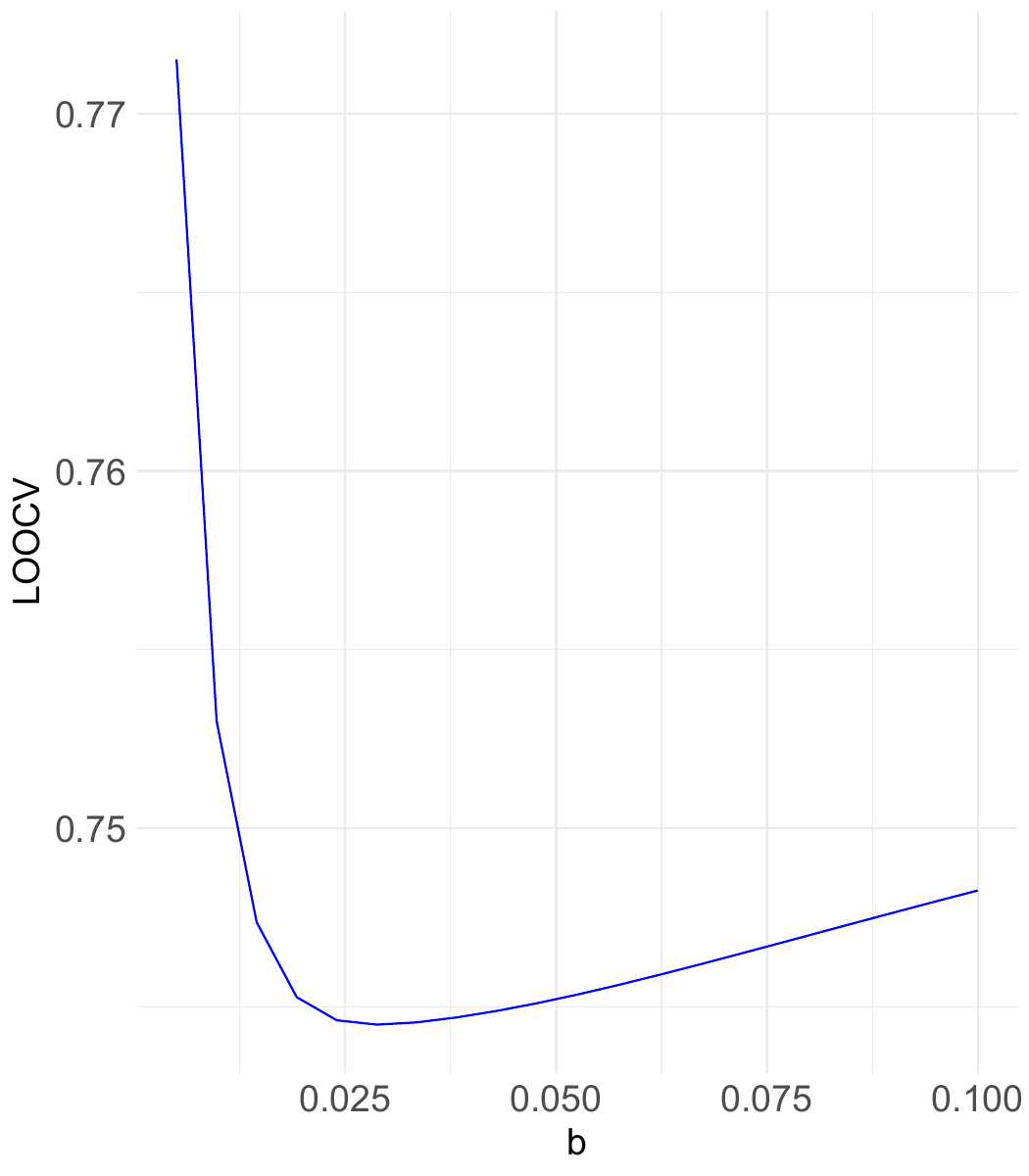}
\includegraphics[width=0.55\textwidth]{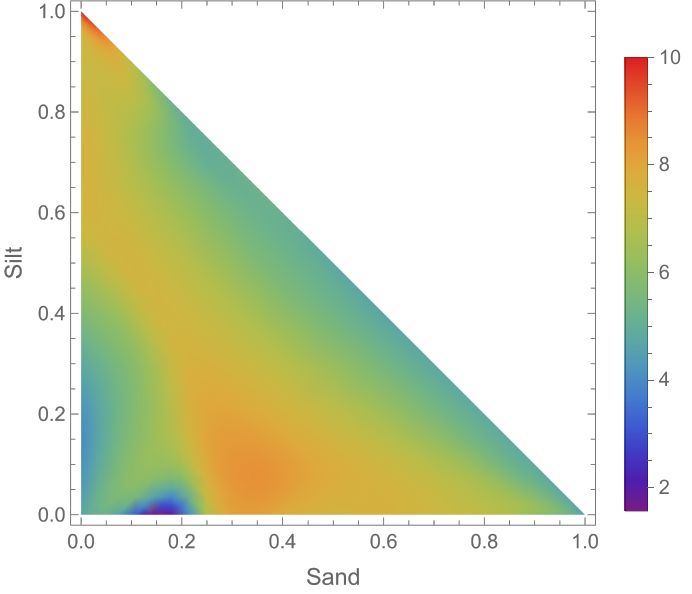}
\caption{Left panel: Plot of the LOOCV criterion against the smoothing parameter for the GEMAS dataset. Right panel: Density plot of the estimated pH in CaCl$_2$ as a function of the proportion of sand and silt.}
\label{fig:density.plot}
\end{figure}

The smoothing parameter $b$ is selected using the leave-one-out cross-validation (LOOCV) criterion described in Section~\ref{sec:simulations}. The left panel of Figure~\ref{fig:density.plot} shows the graph of $\mathrm{LOOCV}$ as a function of the smoothing parameter $b$. The optimal smoothing parameter under this criterion is found to be
\[
\hat{b} = \operatorname*{arg\,min}_{b\in (0,\infty)} \mathrm{LOOCV}(b) \approx 0.0303,
\]
and the resulting density plot of $\smash{\hat{m}_{2083,\hat{b}}^{\mathrm{LL}}}$ is illustrated in the right panel of Figure~\ref{fig:density.plot}.

\newpage
To interpret this plot, it is useful to know the terms typically associated with different ranges of soil pH \citep{Bickelhaupt}:
\begin{multicols}{2}
\begin{itemize}\setlength\itemsep{-0.3em}
\item [a)] Extremely acidic: Less than 4.5
\item [b)] Very strongly acidic: 4.5--5.0
\item [c)] Strongly acidic: 5.1--5.5
\item [d)] Moderately acidic: 5.6--6.0
\item [e)] Slightly acidic: 6.1--6.5
\item [f)] Neutral: 6.6--7.3
\item [g)] Slightly alkaline: 7.4--7.8
\item [h)] Moderately alkaline: 7.9--8.4
\item [i)] Strongly alkaline: 8.5--9.0
\item [j)] Very strongly alkaline: Greater than 9.1
\end{itemize}
\end{multicols}

Based on the density plot and the pH level terms listed above, observations and interpretations regarding pH levels at different sand, silt, and clay compositions are as follows:
\begin{enumerate}\setlength\itemsep{-0.1em}
\item Low sand, very low silt, high clay (e.g., 15--20\% sand, 0--5\% silt, and the balance in clay): This region is in a blueish-purple shade, suggesting a pH level below 4.5, which would classify the soil as \textit{extremely acidic}.

\item Very high sand, very low silt (e.g., 97--100\% sand, 0--5\% silt, and 0--3\% clay): This region is in a reddish-orange shade, suggesting a pH level above 9.1, which would classify the soil as \textit{very strongly alkaline}.

\item High sand, low silt, moderate clay (e.g., 70--80\% sand, 10--20\% silt, and the balance in clay): This region appears mostly yellow-green, which suggests a pH level around 6.5 to 7.5, indicating that the soil is likely \textit{slightly acidic to neutral}.

\item Moderate sand, low silt (e.g., 30--40\% sand, 5--15\% silt, and the balance in clay): This region is in an orange shade, suggesting a pH level around 7.5 to 8.5, indicating that the soil is likely \textit{slightly to moderately alkaline}.

\item High sand and silt, with sand above 20\% (e.g., proportions of sand and silt together exceeding 90\%, and sand alone above 20\%): This region appears in a green shade, suggesting a pH level around 5.6 to 6.5, which would classify the soil as \textit{moderately to slightly acidic}.

\item Very low sand, low to moderate silt (e.g., 0--3\% sand, 5--30\% silt, and the balance in clay): This region is in a blue-green shade, suggesting a pH level around 4.5 to 5.0, indicating that the soil is likely \textit{very strongly acidic}.

\item Moderate sand and silt (e.g., the sum of sand and silt is 60\% $\pm$ 5\%): This region is in a light orange shade, suggesting a pH level around 7.5 to 8.0, which would classify the soil as \textit{slightly to moderately alkaline}.
\end{enumerate}

These interpretations provide a detailed but non-comprehensive view of soil pH trends in relation to different proportions of sand, silt, and clay in Figure~\ref{fig:density.plot}. Other observations and comments could possibly be made by experts in the field.

\begin{appendices}

\renewcommand{\thesection}{Appendix~\Alph{section}}

\section{Proofs of the main results}\label{app:proofs}

\renewcommand{\thesection}{\Alph{section}}

\subsection{Proof of Proposition~\ref{prop:bias}}\label{Appendix:A1}

For each $\bb{s}\in \mathcal{S}_d$, let $\bb{\xi}_{\bb{s}} = (\xi_1, \dots, \xi_d) \sim \mathrm{Dirichlet}\{\bb{s}/b + 1, (1 - \|\bb{s}\|_1)/b + 1\}$, and write
\[
\Bias\{\hat{m}_{n,b}^{\mathrm{GM}}(\bb{s})\}
= \big[\EE\{\hat{m}_{n,b}^{\mathrm{GM}}(\bb{s})\} - \EE\{m(\bb{\xi}_{\bb{s}})\}\big] + \big[\EE\{m(\bb{\xi}_{\bb{s}})\} - m(\bb{s})\big] \equiv \mbox{(I)} + \mbox{(II)}.
\]
By the Lipschitz continuity of $m$ due to Assumption~\ref{ass:1}, together with Assumptions~\ref{ass:3}, \ref{ass:4.a}, \ref{ass:4.b}, \ref{ass:4.c}, and \ref{ass:4.d}, one easily finds that, as $n \to \infty$,
\[
|\mbox{(I)}| \leq \sum_{i=1}^n \int_{B_i} |m(\bb{x}_i) - m(\bb{x})| \kappa_{\bb{s},b}(\bb{x}) \, \rd \bb{x} \ll n^{-1/d} \sum_{i=1}^n \int_{B_i} \kappa_{\bb{s},b}(\bb{x}) \, \rd \bb{x} \ll n^{-1/d} = \oo(b).
\]
Also, given that $m$ is twice continuously differentiable by Assumption~\ref{ass:1}, one deduces that $\mbox{(II)} = b g(\bb{s}) + \oo(b)$ as $n \to \infty$ \cite[Theorem~1]{MR4319409}, where $g(\bb{s})$ is specified in \eqref{eq:fonction-g}. It suffices to combine the asymptotic behavior of (I) and (II) to complete the proof.

\subsection{Proof of Proposition~\ref{prop:var}}\label{Appendix:A2}

Fix a set $\mathcal{J} \subseteq [d]$ and $\bb{\lambda}\in (2, \infty)^d$. Assume that $\bb{s} = \bb{s}(b)$ is chosen inside $\mathrm{Int}(\mathcal{S}_d)$ in such a way that $s_j / b \to \lambda_j$ for all $j \in \mathcal{J}$, and $s_j$ is fixed for all $j \in [d] \backslash \mathcal{J}$. For each integer $i \in [n]$, the multivariate mean value theorem ensures the existence of a vector $\bb{t}_i\in B_i$ such that $\smash{\int_{B_i} \kappa_{\bb{s},b}(\bb{t}) \, \rd \bb{t} = \lambda(B_i) \kappa_{\bb{s},b}(\bb{t}_i)}$. Given that the $Y_i$'s are independent with $\Var(Y_i)=\sigma^2(\bb{x}_i)$ for each $i\in [n]$, the Lipschitz continuity of $\sigma^2$ from Assumption~\ref{ass:1}, together with Assumption~\ref{ass:4.d}, implies that
\begin{equation}\label{eq:variance.expansion}
\Var\{\hat{m}^{\mathrm{GM}}_{n,b}(\bb{s})\}
= \sum_{i=1}^n \sigma^2(\bb{x}_i)\left\{\int_{B_i}\kappa_{\bb{s},b}(\bb{x}) \, \rd \bb{x}\right\}^2
= \sum_{i=1}^n \lambda(B_i)\{\sigma^2(\bb{t}_i) + \OO(n^{-1/d})\} \lambda(B_i) \kappa_{\bb{s},b}^2(\bb{t}_i).
\end{equation}
The positivity and Lipschitz continuity of $1/f$ from Assumptions~\ref{ass:1} and~\ref{ass:2}, together with Assumptions~\ref{ass:4.d} and~\ref{ass:4.e}, ensure that, for each integer $i\in[n]$, there exists $\bb{\xi}_i\in B_i$ (here, the multivariate mean value theorem is used again) such that
\[
\lambda(B_i) = \frac{1}{n f(\bb{\xi}_i)} + \OO(n^{-1-1/d}) = \frac{1}{n f(\bb{t}_i)} + \OO(n^{-1-1/d}).
\]
Then, Eq.~\eqref{eq:variance.expansion}, the positivity and Lipschitz continuity of $\sigma^2/f$ from Assumptions~\ref{ass:1} and~\ref{ass:2}, as well as Lemma~\ref{lem:weighted.A.b}, enable one to show that
\begin{equation}\label{eq:var.expansion}
\begin{aligned}
\Var\{\hat{m}^{\mathrm{GM}}_{n,b}(\bb{s})\}
&= \sum_{i=1}^n \lambda(B_i)\frac{\sigma^2(\bb{t}_i)}{nf(\bb{t}_i)}\kappa_{\bb{s},b}^2(\bb{t}_i)
+ \OO(n^{-1-1/d})\sum_{i=1}^n \lambda(B_i)\kappa_{\bb{s},b}^2(\bb{t}_i) \\
&= \frac{\sigma^2(\bb{s})}{nf(\bb{s})}A_b(\bb{s})\{1 + \OO(b^{1/2})\} \{1 + \oo(1)\} + \OO(n^{-1-1/d})A_b(\bb{s}) \{1 + \oo(1)\}.
\end{aligned}
\end{equation}
The conclusion follows from Lemma~\ref{lem:A.b.asymptotics}.

\subsection{Proof of Theorem~\ref{thm:MISE}}

Given that $m$ is twice continuously differentiable on $\mathcal{S}_d$ by Assumption~\ref{ass:1} and $\mathcal{S}_d$ is compact, the second-order partial derivatives of $m$ are bounded. Therefore, using Proposition~\ref{prop:bias} together with Lebesgue's dominated convergence theorem, one finds that
\[
b^{-2} \int_{\mathcal{S}_d} \Bias\{\hat{m}^{\mathrm{GM}}_{n,b}(\bb{s})\}^2 \, \rd \bb{s}
= \int_{\mathcal{S}_d} g^2(\bb{s}) \, \rd \bb{s} + \oo(1).
\]
Similarly, using the upper bound and almost-everywhere convergence of $b^{d/2} A_b$ from Lemma~\ref{lem:A.b.asymptotics}, the positivity and Lipschitz continuity of $\sigma^2/f$ from Assumptions~\ref{ass:1} and~\ref{ass:2}, and Lebesgue's dominated convergence theorem, one deduces from Eq.~\eqref{eq:var.expansion} in the proof of Proposition~\ref{prop:var} that
\[
n b^{d/2} \int_{\mathcal{S}_d} \Var\big\{\hat{m}^{\mathrm{GM}}_{n,b}(\bb{s})\big\} \, \rd \bb{s}
= \int_{\mathcal{S}_d} \frac{b^{d/2} A_b(\bb{s}) \sigma^2(\bb{s})}{f(\bb{s})} \, \{1 + \oo(1)\} \, \rd \bb{s}
= \int_{\mathcal{S}_d} v(\bb{s}) \, \rd \bb{s} + \oo(1).
\]
Combining the last two equations makes it possible to conclude.

\subsection{Proof of Theorem~\ref{thm:CLT}}

Let $\bb{s}\in \mathrm{Int}(\mathcal{S}_d)$ be given, and let $R_n = n^{1/2} b^{d/4}$. Note that Assumption~\ref{ass:3} implies $R_n\to \infty$ as $n \to \infty$. Moreover, simple algebraic manipulations show that the restriction $d\in \{1,2,3\}$ is necessary to satisfy Assumption~\ref{ass:3} under the conditions $R_n b \to 0$ and $n^{2/(d+4)} \, b \to \nu$ in $(i)$~and~$(ii)$, respectively. Now, consider the decomposition
\[
\{\hat{m}_{n,b}^{\mathrm{GM}}(\bb{s}) - m(\bb{s})\} = \big[\hat{m}_{n,b}^{\mathrm{GM}}(\bb{s}) - \EE\{\hat{m}_{n,b}^{\mathrm{GM}}(\bb{s})\}\big] + \big[\EE\{\hat{m}_{n,b}^{\mathrm{GM}}(\bb{s})\} - m(\bb{s})\big].
\]
By Proposition~\ref{prop:bias}, the rescaled bias term $\smash{n^{1/2} b^{d/4} [\EE\{\hat{m}_{n,b}^{\mathrm{GM}}(\bb{s})\} - m(\bb{s})]}$ converges to $0$ if $R_n b\to 0$, and $\smash{n^{2/(d+4)} [\EE\{\hat{m}_{n,b}^{\mathrm{GM}}(\bb{s})\} - m(\bb{s})]}$ converges to $\nu g(\bb{s})$ if $n^{2/(d+4)} b\to \nu$. Therefore, to conclude the proof, it suffices to show that, as $n \to \infty$,
\begin{equation}\label{eq:to.prove}
R_n \big[\hat{m}_{n,b}^{\mathrm{GM}}(\bb{s}) - \EE\{\hat{m}_{n,b}^{\mathrm{GM}}(\bb{s})\}\big] \rightsquigarrow \mathcal{N} \big[ 0, v(\bb{s}) \big],
\end{equation}
because $(i)$ and $(ii)$ follow immediately by Slutsky's lemma (and a rescaling of the variance for $(ii)$).

Recalling the definition of the D-GM estimator in \eqref{eq:GM.estimator}, and noting that the errors $\e_1,\dots,\e_n$ have mean zero and are mutually independent, one has the representation
\[
\hat{m}_{n,b}^{\mathrm{GM}}(\bb{s}) - \EE\{\hat{m}_{n,b}^{\mathrm{GM}}(\bb{s})\} = \sum_{i=1}^n Z_{\bb{s},b,i}, \qquad \text{with } Z_{\bb{s},b,i} = \e_i \int_{B_i} \kappa_{\bb{s},b}(\bb{x}) \, \rd \bb{x},
\]
where the random variables $Z_{\bb{s},b,1},\dots,Z_{\bb{s},b,n}$ have mean zero and are mutually independent. The asymptotic normality in Eq.~\eqref{eq:to.prove} will be established by verifying the following Lindeberg condition for triangular arrays \citep[see, e.g.,][p.~359]{MR1324786}:
\begin{equation}\label{eq:Lindeberg}
\forall_{\eta\in (0,\infty)} \qquad \lim_{n\to \infty} \frac{1}{s_{n,b}^2} \sum_{i=1}^n \EE\big(|Z_{\bb{s},b,i}|^2 \, \ind_{\{|Z_{\bb{s},b,i}| > \eta \, s_{n,b}\}}\big) = 0,
\end{equation}
where $s_{n,b}^2 = \sum_{i=1}^n \EE(|Z_{\bb{s},b,i}|^2)$.

Using Proposition~\ref{prop:var} in the case where $\bb{s}\in \mathrm{Int}(\mathcal{S}_d)$ is fixed (i.e., $\mathcal{J} = \emptyset$), one has
\[
s_{n,b}^2 = R_n^{-2} v(\bb{s}) \{1 + \oo_{\bb{s}}(1)\}.
\]
Additionally, using Assumptions~\ref{ass:2}~and~\ref{ass:4.e} to obtain $\lambda(B_i) \ll n^{-1}$, and calling on the upper bound on the Dirichlet kernel from Lemma~\ref{lem:local.bound} to get $\max\{\kappa_{\bb{s},b}(\bb{x}) : \bb{x}\in \mathcal{S}_d\} \ll b^{-d/2} \psi_{\emptyset}(\bb{s})$, one finds, for every integer $i\in [n]$,
\[
|Z_{\bb{s},b,i}| \leq |\e_i| \, \lambda(B_i) \max\{\kappa_{\bb{s},b}(\bb{x}) : \bb{x}\in \mathcal{S}_d\} \ll |\e_i| \, R_n^{-2} \psi_{\emptyset}(\bb{s}).
\]
Therefore, combining the last two equations, it suffices to show that
\[
\forall_{\delta\in (0,\infty)} \qquad \lim_{n \to \infty} \frac{1}{R_n^{-2}} \sum_{i=1}^n R_n^{-4} \, \EE\big(\e_i^2 \, \ind_{\raisebox{-0.1mm}[0pt][0pt]{$\scriptstyle \{|\e_i| \, R_n^{-2} > \delta R_n^{-1}\}$}}\big) = 0,
\]
which is equivalent to the assumption~\eqref{eq:asympotic.normality.condition}. This establishes Eq.~\eqref{eq:Lindeberg}, which in turn implies Eq.~\eqref{eq:to.prove} and completes the proof.

\renewcommand{\thesection}{Appendix~\Alph{section}}

\section{Technical lemmas}\label{app:tech.lemmas}

\renewcommand{\thesection}{\Alph{section}}

The first lemma gives a uniform upper bound on the Dirichlet kernel.

\begin{lemma}[{\cite[Lemma~2]{MR4319409}}]\label{lem:local.bound}
If $\alpha_1, \dots, \alpha_d, \beta\in [2,\infty)$, then
\[
\sup\{K_{\bb{\alpha},\beta}(\bb{x}) : \bb{x}\in \mathcal{S}_d\} \leq \sqrt{\frac{\|\bb{\alpha}\|_1 + \beta - 1}{(\beta - 1) \prod_{j\in [d]} (\alpha_j - 1)}} ~ (\|\bb{\alpha}\|_1 + \beta - d - 1)^d.
\]
\end{lemma}

The second lemma provides the asymptotic behavior of $A_b(\bb{s}) = \int_{\mathcal{S}_d}\kappa_{\bb{s},b}^2(\bb{t}) \, \rd \bb{t}$ as $b \to 0$, whether the sequence of points $\bb{s} = \bb{s}(b)$ lies well inside or near the boundary of $\mathcal{S}_d$.

\begin{lemma}[{\cite[Lemma~1]{MR4319409}}]\label{lem:A.b.asymptotics}
As $b \to 0$, one has, uniformly in $\bb{s}\in \mathcal{S}_d$,
\[
0 < A_b(\bb{s}) \leq \frac{b^{(d + 1) / 2} \, (1 / b + d)^{d + 1/2}}{(4\pi)^{d/2} \sqrt{(1 - \|\bb{s}\|_1) \prod_{j=1}^d s_j}} \, \{1 + \OO(b)\}.
\]
For any vector $\bb{\lambda} = (\lambda_1,\dots,\lambda_d)\in (0,\infty)^d$ and subset $\mathcal{J}\subseteq [d]$ of indices, one also has, as $b\to 0$,
\[
A_b(\bb{s}) = b^{-d/2} \, \psi_{\emptyset}(\bb{s}) \times \{1 + \OO_{\bb{s}}(b) \}
\]
if $s_j$ is fixed for every integer $j \in [d]$, while
\[
A_b (\bb{s}) = b^{-(d + |\mathcal{J}|)/2} \psi_{\mathcal{J}}(\bb{s}) \prod_{j \in \mathcal{J}} \frac{\Gamma(2\lambda_j + 1)}{2^{2\lambda_j + 1} \Gamma^2(\lambda_j + 1)} \, \{1 + \OO_{\bb{\lambda},\bb{s}}(b)\},
\]
if $s_j / b \to \lambda_j$ for every integer $j \in \mathcal{J}$, and $s_j$ is fixed for all $j \in [d] \backslash \mathcal{J}$.
\end{lemma}

Next, Lemma~\ref{lem:A.b.asymptotics} is extended to the asymptotics of $\int_{\mathcal{S}_d} G(\bb{t}) \kappa_{\bb{s},b}^2(\bb{t}) \, \rd \bb{t}$ for a Lipschitz weight $G$.

\begin{lemma}\label{lem:weighted.A.b}
Let $\mathcal{J} \subseteq [d]$ and $\bb{\lambda} = (\lambda_1,\dots,\lambda_d)\in (2,\infty)^d$ be given. Let $\bb{s} = \bb{s}(b)$ be a sequence selected inside $\mathrm{Int}(\mathcal{S}_d)$ such that $s_j / b \to \lambda_j$ for all $j \in \mathcal{J}$, and $s_j$ is fixed for all $j \in [d]\backslash \mathcal{J}$. Suppose that the function $G$ is positive and Lipschitz on $\mathcal{S}_d$.
\begin{itemize}\setlength\itemsep{0em}
\item[(i)] If $b = b(n)\to 0$ as $n\to\infty$, then $\int_{\mathcal{S}_d} G(\bb{t})\kappa_{\bb{s},b}^2(\bb{t}) \, \rd \bb{t} = G(\bb{s})A_b(\bb{s})\{1 + \OO(b^{1/2})\}$.
\item[(ii)] If Assumptions~\ref{ass:3},~\ref{ass:4.a},~\ref{ass:4.c},~\ref{ass:4.d} hold, and if, for every integer $i\in [n]$, $\bb{\tau}_i\in B_i$ is such that $\int_{B_i} \kappa_{\bb{s},b}(\bb{t}) \, \rd \bb{t} = \lambda(B_i) \kappa_{\bb{s},b}(\bb{\tau}_i)$, then
\[
\sum_{i=1}^n \lambda(B_i)G(\bb{\tau}_i)\kappa_{\bb{s},b}^2(\bb{\tau}_i)
= \int_{\mathcal{S}_d} G(\bb{t})\kappa_{\bb{s},b}^2(\bb{t}) \, \rd \bb{t} \times \{1 + \oo(1)\} \qquad \text{as } n\to\infty.
\]
\end{itemize}
\end{lemma}

\begin{proof}[Proof]
For the proof of $(i)$, see \citet[Eq.~(15)]{MR4319409}, where what is denoted $G$ here is denoted $f$ there. To prove $(ii)$, let $L_G\in (0,\infty)$ denote a Lipschitz constant of $G$ on $\mathcal{S}_d$, and set $M_G = \max \{|G(\bb{t})| : \bb{t}\in \mathcal{S}_d\} < \infty$. By the assumption on $\bb{\tau}_i$ and the triangle inequality,
\[
\begin{aligned}
\left|\sum_{i=1}^n \lambda(B_i)G(\bb{\tau}_i)\kappa_{\bb{s},b}^2(\bb{\tau}_i) - \sum_{i=1}^n G(\bb{\tau}_i) \int_{B_i} \kappa_{\bb{s},b}^2(\bb{t}) \, \rd \bb{t}\right|
&= \left|\sum_{i=1}^n G(\bb{\tau}_i) \int_{B_i} \{\kappa_{\bb{s},b}(\bb{\tau}_i) - \kappa_{\bb{s},b}(\bb{t})\} \kappa_{\bb{s},b}(\bb{t}) \, \rd \bb{t}\right| \\
&\leq M_G \sum_{i=1}^n \int_{B_i} |\kappa_{\bb{s},b}(\bb{\tau}_i) - \kappa_{\bb{s},b}(\bb{t})| \kappa_{\bb{s},b}(\bb{t}) \, \rd \bb{t}.
\end{aligned}
\]
From the multivariate mean value theorem, for every integer $i\in [n]$ and every $\bb{t}\in B_i$, there exists $\bb{\eta}_{i,\bb{t}}$ on the line segment joining $\bb{\tau}_i$ and $\bb{t}$ (hence $\bb{\eta}_{i,\bb{t}}\in \mathcal{S}_d$, because $\mathcal{S}_d$ is convex) such that
\[
\left|\kappa_{\bb{s},b}(\bb{\tau}_i) - \kappa_{\bb{s},b}(\bb{t})\right| \leq \sum_{j=1}^d \left|\Big.\frac{\partial}{\partial r_j} \, \kappa_{\bb{s},b}(\bb{r})\Big|_{\bb{r} = \bb{\eta}_{i,\bb{t}}}\right| |\tau_{i,j} - t_j |,
\]
where $\tau_{i, j}$ and $t_j$ denote the $j$th components of $\bb{\tau}_i$ and $\bb{t}$, respectively. For $j \in [d]$, let $\bb{e}_j$ denote the $j$th standard basis vector in $\R^d$. Recall that, for any reals $\alpha_1, \dots, \alpha_d, \beta\in (1, \infty)$ and $\bb{x}\in \mathcal{S}_d$, one has
\[
\frac{\partial}{\partial x_j} \, K_{\bb{\alpha},\beta}(\bb{x}) = (\|\bb{\alpha}\|_1+\beta-1) \{K_{\bb{\alpha} - \bb{e}_j,\beta}(\bb{x}) - K_{\bb{\alpha},\beta-1}(\bb{x})\}.
\]
Choosing $\bb{\alpha} = \bb{s}/b + \bb{1}$ and $\beta = (1 - \|\bb{s}\|_1)/b + 1$, and applying the upper bound on the Dirichlet kernel from Lemma~\ref{lem:local.bound} (note that $\bb{s}/b + \bb{1} - \bb{e}_j \in [2,\infty)^d$ and $(1 - \|\bb{s}\|_1)/b\in [2,\infty)$ for $b$ small enough), one has
\[
\begin{aligned}
\max_{\bb{r}\in \mathcal{S}_d} \left|\frac{\partial}{\partial r_j} \, \kappa_{\bb{s},b}(\bb{r})\right|
&= \left(\frac{1}{b} + d\right) \, \max_{\bb{r} \in \mathcal{S}_d} \left|K_{\bb{s}/b + \bb{1} - \bb{e}_j, (1 - \|\bb{s}\|_1)/b + 1}(\bb{r}) - K_{\bb{s}/b + \bb{1}, (1 - \|\bb{s}\|_1)/b}(\bb{r})\right| \\
&\ll b^{-1} \times b^{-(d + |\mathcal{J}|)/2} \psi_{\mathcal{J}}(\bb{s}).
\end{aligned}
\]
Given that $\max_{j \in [d]} |\tau_{i, j} - t_j| \ll n^{-1/d}$ by Assumption~\ref{ass:4.d}, it follows from the previous displayed equations that
\[
\begin{aligned}
\left|\sum_{i=1}^n \lambda(B_i)G(\bb{\tau}_i)\kappa_{\bb{s},b}^2(\bb{\tau}_i) - \sum_{i=1}^n G(\bb{\tau}_i) \int_{B_i} \kappa_{\bb{s},b}^2(\bb{t}) \, \rd \bb{t}\right|
&\ll b^{-(d + |\mathcal{J}|)/2} \psi_{\mathcal{J}}(\bb{s}) \times b^{-1} n^{-1/d} \times \sum_{i=1}^n \int_{B_i} \kappa_{\bb{s},b}(\bb{t}) \, \rd \bb{t} \\
&\ll b^{-(d + |\mathcal{J}|)/2} \psi_{\mathcal{J}}(\bb{s}) \times b^{-1} n^{-1/d} \\[2mm]
&= \oo\{A_b(\bb{s})\},
\end{aligned}
\]
where the equality follows from Assumption~\ref{ass:3} and Lemma~\ref{lem:A.b.asymptotics}.

Moreover, by the Lipschitz continuity of $G$ and Assumption~\ref{ass:4.d},
\[
\begin{aligned}
\left|\sum_{i=1}^n G(\bb{\tau}_i) \int_{B_i} \kappa_{\bb{s},b}^2(\bb{t}) \, \rd \bb{t} - \int_{\mathcal{S}_d} G(\bb{t}) \kappa_{\bb{s},b}^2(\bb{t}) \, \rd \bb{t}\right|
&\leq \sum_{i=1}^n \int_{B_i} |G(\bb{\tau}_i) - G(\bb{t})| \kappa_{\bb{s},b}^2(\bb{t}) \, \rd \bb{t} \\
&\leq L_G \max_{i\in [n]} \, \mathrm{diam}(B_i) \int_{\mathcal{S}_d} \kappa_{\bb{s},b}^2(\bb{t}) \, \rd \bb{t}
\ll \frac{A_b(\bb{s})}{n^{1/d}} = \oo\{A_b(\bb{s})\}.
\end{aligned}
\]
Combining the last two displayed equations proves $(ii)$ given that $G$, being positive, takes values in a compact subset of $(0,\infty)$. This concludes the proof.
\end{proof}

\end{appendices}

\section*{Reproducibility}\label{sec:reproducibility}
\addcontentsline{toc}{section}{Reproducibility}

The \textsf{R} code that generated the figures, the simulation results and the data analysis are available online in the GitHub repository of \citet{DaayebGenestKhardaniKlutchnikoffOuimet2026github}.

\section*{Acknowledgments}
\addcontentsline{toc}{section}{Acknowledgments}

The simulations in Section~\ref{sec:simulations} were carried out using the computational resources supplied by Calcul Qu\'ebec (\href{https://www.calculquebec.ca}{www.calculquebec.ca}) and the Digital Research Alliance of Canada (\href{https://www.alliancecan.ca}{www.alliancecan.ca}). The authors thank the reviewers for their constructive comments, which have significantly enhanced the quality of this paper. Thanks are also due to Raimon Tolosana-Delgado (Helmholtz-Zentrum Dresden-Rossendorf, Germany) for providing the GEMAS dataset used for illustration in Section~\ref{sec:application}.

\section*{Funding}
\addcontentsline{toc}{section}{Funding}

Genest's work was funded by the Natural Sciences and Engineering Research Council of Canada (NSERC) Grant RGPIN-2024-04088 and the Canada Research Chairs Program (Grant 950-231937). Ouimet acknowledges funding from NSERC through Discovery Grant RGPIN-2026-04471 and Discovery Launch Supplement DGECR-2026-00449.

\phantomsection
\addcontentsline{toc}{section}{References}

\setlength{\bibsep}{0pt plus 0ex}

\bibliographystyle{apalike}
\bibliography{bib}

\end{document}